\documentclass[leqno]{amsart}

\usepackage{latexsym,amssymb,amsthm,amsmath,amscd,xcolor}

\theoremstyle{plain}
\newtheorem{theorem}{Theorem}[section]

\newtheorem{proposition}{Proposition}[section]

\newtheorem{definition}{Definition}[section]

\numberwithin{equation}{section}

\theoremstyle{remark}
\newtheorem{remark}{Remark}[section]

 \numberwithin{equation}{section}

\newtheorem*{Theorem A}{{\bf Theorem A}}
\newtheorem*{Theorem B}{{\bf Theorem B}}
\newtheorem*{Theorem C}{Theorem C}

 \numberwithin{equation}{section}
\def\<{\left < }
\def\>{\right >}
\def\({\left ( }
\def\){\right )}
\def\e{\eqref}

\begin{document}

\title[$n$-harmonicity, minimality, conformality  and cohomology] {$n$-harmonicity, minimality, conformality  and cohomology}

\author[B.-Y. Chen]{Bang-Yen Chen}
\address{Department of Mathematics\\
    Michigan State University \\East Lansing, Michigan 48824--1027\\ U.S.A.}
\email{chenb@msu.edu}

\author[S. W. Wei]{Shihshu Walter Wei}
\address{Department of Mathematics\\
University of Oklahoma\\ Norman, Oklahoma 73019-0315\\ U.S.A.}
\email{wwei@ou.edu}

\begin{abstract}  By studying cohomology classes that are related with $n$-harmonic morphisms and $F$-harmonic maps, we augment and extend several results on $F$-harmonic maps, harmonic maps in \cite {A, BE, EL2}, $p$-harmonic morphisms in \cite {OW}, and also revisit our previous results in \cite{CW08,CW09,W9}  on Riemannian submersions and $n$-harmonic morphisms which are submersions. The results, for example Theorem \ref{T:3.2} obtained by utilizing the $n$-conservation law \eqref{2.6}, are sharp.
\end{abstract}





 
\thanks{}

\date{}

\maketitle




\section{Introduction}\label{S1}

Harmonicity and its variants are related with the topology and geometry of manifolds. It was shown in \cite {S.W.Wei 1} that homotopy classes can be represented by $p$-harmonic maps
 (see, e.g. \cite {W9}, for definition and examples of $p$-harmonic maps):
\vskip.1in

\noindent {\bf Theorem A.} {\it If $N^{n}$ is a compact Riemannian $n$-manifold, then for any positive integer $i$, each class in the i-th homotopy group $\pi_i(N^{n})$ can be represented by a $C^{1,\alpha}$ $p$-harmonic map $u_0$ from an $i$-dimensional sphere $S^i$ into $N^{n}$ minimizing $p$-energy in its homotopy class for any $p>i$.}
\vskip.1in


On the other hand, B.-Y. Chen established in \cite{c2005} the following result involving Riemannian submersion, minimal immersion, and  cohomology class.
\vskip.1in

\noindent {\bf Theorem B.} (\cite{c2005}) {\it Let  $\pi : (M^{m},g_{M})\to
(B^{b},g_{B})$ be a Riemannian submersion with minimal fibers and orientable
base manifold $B^{b}$. If  $M^{m}$ is a closed manifold with cohomology
class $H^{b}(M^{m},\mathbf R) =0$, then the horizontal distribution
$\mathcal{H}$ of the Riemannian submersion is never integrable.
Thus the submersion $\pi$ is never non-trivial.}
\vskip.1in

Whereas $p$-harmonic maps represent homotopy classes,  B.-Y. Chen and S.\,W. Wei connected  the two seemingly unrelated areas of $p$-harmonic morphisms and cohomology classes in the following.
\vskip.1in

\noindent {\bf Theorem C.} (\cite {CW08,CW09}) {\it Let $u: (M^{m},g_{M}) \to (N^{n}, g_{N})$ be an $n$-harmonic morphism which is a submersion. If $N^{n}$ is an orientable $n$-manifold and $M^{m}$ is a closed $m$-manifold with $n$-th cohomology class $H^n(M^{m},\mathbf R) = 0$, then the horizontal distribution $\mathcal H$ of $u$ is never integrable. Hence the submersion $u$ is always non-trivial.}
\vskip.1in

This recaptures Theorem B when $\pi : M^{m}\to B^{b}$ is  a  Riemannian submersion with minimal fibers
and orientable base manifold $B^{b}$. While a horizontally weak conformal $p$-harmonic map is a
$p$-harmonic morphism (cf. e.g., \cite [Theorem 4]{CW09}), $p$-harmonic morphism is also linked to cohomology class as follows.
\vskip.1in

\noindent {\bf Theorem D.} (\cite {CW08,CW09}) {\it Let $u: (M^{m},g_{M}) \to (N^{n},g_{N})$ be an $n$-harmonic morphism which is a submersion. Then the pull back of the volume element of the base manifold $N^{n}$ is a harmonic $n$-form if and only if the horizontal distribution $\mathcal H$ of $u$ is completely integrable.}
\vskip.1in

Following the proofs given in \cite {CW08,CW09}, and by applying a characterization theorem of a $p$-harmonic morphism from \cite{BG,BL}, and \cite [Theorem 2.5] {W9}, we seek a dual version of Theorem D. In particular,  $p$-harmonic maps and cohomology classes are interrelated in \cite {W9} as follows.
\vskip.1in

\noindent {\bf Theorem E.} {\it  Let $M^{m}$ be a closed $m$-manifold and  $u: (M^{m},g_{M}) \to (N^{n},g_{N})$ be an $n$-harmonic map which is a submersion. If $M^{m}$ is a closed $m$-manifold and the horizontal distribution $\mathcal H$ of $u$ is integrable and $u$ is an $n$-harmonic morphism, then we have $H^n(M,\mathbf R) \ne 0$.}
\vskip.1in

\noindent {\bf Theorem F.} (\cite {W9}) {\it Let $u: (M^{m},g_{M}) \to (N^{n},g_{N})$ be an $n$-harmonic map  which is a submersion such that the horizontal distribution $\mathcal H$ of $u$ is integrable. If $M^{m}$ is a closed manifold with cohomology class $H^n(M^{m},\mathbf R) = 0$. Then $u$ is not an $n$-harmonic morphism. Thus the submersion $u$ is always nontrivial.}
\vskip.1in

The purpose of this paper is to point out the underlying essence of the foregoing Theorems C, D, E, and F is an application of stress-energy tensor and a conservation law.
The results, for example Theorem \ref{T:3.2} obtained by utilizing the $n$-conservation law \eqref{2.6}, are sharp.

\section{Preliminaries}\label{S2}

\subsection{Submersions}\label{S2.1}

A differential map $u: (M^{m},g_{M}) \to (N^{n},g_{N})$  between two Riemannian manifolds  is called a submersion at a point $x\in M^{m}$ if its differential $du_{x}: T_x(M^{m})\to T_{u(x)}(N^{n})$ is a surjective linear map. A differentiable map $u$ that is a submersion at each point $x\in M^{m}$ is called a {\it submersion}. For each point $x\in N^{n}$, $u^{-1}(x)$ is called a {\it fiber}.
For a submersion $u: M \to N$, let $\mathcal H_x$ denote the orthogonal complement of Kernel$\, \big (du_x : T_x(M^{m})\to T_{u(x)}(N^{n})\big )$. Let $\mathcal H=\{\mathcal H_x: x\in M^{m}\}$ denote the horizontal distribution of $u$. 

A submersion $u: (M^{m},g_{M}) \to (N^{n},g_{N})$ is called {\it horizontally weakly conformal} if the restriction of $du_x$ to $\mathcal H _x$ is conformal, i.e., there exists a smooth function $\lambda$ on $M^{m}$ such that 
\begin{align}\label{2.1} u^{*} g_{N} = \lambda^2  {g _M}_ {|_{\mathcal H}} \;\;\;\; {\rm or}\;\; \;\;
g_{N}\big (du_x(X),du_x(Y)\big) = \lambda ^2 (x) g_{M}(X, Y)\end{align}
 for all $X,Y \in \mathcal{H} _x$ and $x \in M^{n}$. 
  If the function $\lambda$ in \eqref{2.1} is  positive, then $u$ is called  {\it horizontally conformal}
and   $\lambda$ is called the {\it dilation} of $u$.  For a horizontally conformal submersion $u$ with dilation $\lambda$, the {\it energy density} of $u$ is  $e_u = \frac{1}{2} {n\lambda ^2}$ (cf. \eqref{2.2}). 

A horizontally conformal submersion with dilation $\lambda \equiv 1$ is called a {\it Riemannian submersion}. 
Recall that a $k$-form $\omega$ on a compact Riemannian manifold is called {\it harmonic} if $\omega$ is both closed and co-closed, i.e., $d\omega=\delta \omega=0$.

In \cite {W9, BG}, generalizing the work of P. Baird and J. Eells for the case $n=2$, and the necessary condition for the fibers being minimal,  
S. W. Wei linked $p$-harmonicity for every $p>1$,  and P. Baird and S. Gudmundsson linked $n$-harmonicity, $n=p=\dim N$ with minimal fibers as follows.

\begin{theorem} [\cite {W9}, Theorem 2.5 ]   Let $u : M
\to N$ be a Riemannian submersion. Then $u$ is a $p$-harmonic map, for every $p
> 1$, if and only if  all fibers $u^{-1}(y)\, ,$ $y\in N$ are $\emph {minimal}$ submanifolds in $ M\, . $\label{T: 2.1}\end{theorem}

\begin{proposition} [\cite {W9}, Proposition 2.4 ]  Let $u : M
\to N$ be a Riemannian submersion. Then $u$ is a $p$-harmonic morphism, for every $p
> 1$, if and only if  all fibers $u^{-1}(y)\, ,$ $y\in N$ are $\emph {minimal}$ submanifolds in $ M\, . $ \label{P: 2.1}
\end{proposition}
The case $p=2$ in Theorem \ref{T: 2.1} and Proposition \ref{P: 2.1} are due to Eells-Sampson \cite{ES}. 

\begin{theorem} [P. Baird and S. Gudmundsson \cite {BG}, Corollary 2.6 ]   If $u: (M^{m},g_{M}) \to (N^{n},g_{N})$ is a horizontally conformal submersion from a Riemannian manifold $M^{m}$ onto a Riemannian manifold $N^{n}$, then $u$ is $n$-harmonic if and only if the fibers of $u$ are minimal in $M^{m}$. \label{T:2.1}\end{theorem}

\begin{remark} (i). The results of linking $p$-harmonicity for every $p>1$,  with minimal fibers in Theorem \ref{T: 2.1} can be extended to $p=1=n$ with minimal fibers. We refer 
to the celebrated work of E. Bombieri - E. De Gorgi - E. Jiusti on minimal cones and the Bernstein problem (\cite {BDG}), S.W. Wei on $1$-harmonic functions (\cite{W7}), P. Baird - S. Gudmundsson on $p$-harmonic maps and minimal submanifods (\cite{BG}), Y.I. Lee - S.W. Wei - A.N. Wang on a generalized $1$-harmonic equation and the inverse mean curvature flow (\cite{LWW}), etc.\quad (ii). We also note that utilizing symmetry, Wu-Yi Hsiang pioneered the study of the inverse image of minimal submanifolds being minimal under appropriate conditions (\cite {Hs}), which marked the birth of {\it equivariant differential geometry} (cf. e.g. W.Y. Hsiang - H.B. Lawson \cite {HL}, S.W. Wei \cite {W3}, etc.). \label{R: 2.1}
\end{remark}
  
\subsection{$F$- and $p$-harmonic morphisms}\label{S2.2}

Let $u: (M^{m},g_{M}) \to (N^{n},g_{N})$ be a differential map between two Riemannian manifolds $M$ and $N$. Denote  $e_u$ the {\it energy density} of  $u$,  which is given by 
\begin{equation}\label{2.2} e_u = \frac 12\sum _{i=1}^m g_{N}\big (du(e_i),du(e_i)\big ) = \frac 12 |du|^2\, , \end{equation}
where $\{e_1, \cdots, e_m\}$ is a local orthonormal frame field on $M^{m}$ and $|du|$ is the Hilbert-Schmidt norm of $du$, determined by the metric $g_{M}$ of $M$ and the metric $g_{N}$ of $N^{n}$. 
{\it The energy of $u$}, denoted by $E(u)$,  is defined to be $$E(u) = \int_M e_u\, dv_g.$$ A smooth map $u: M^{m} \to N^{n}$ is called {\it harmonic} if $u$ is a critical point of the energy functional $E$
with respect to any compactly supported variation.  

Let $F : [0,\infty) \to [0,\infty)$ be a strictly increasing function with  $F(0) = 0$
and let  $u: (M,g_{M}) \to (N,g_{N})$ be a smooth map between two compact Riemannian manifolds. Then the map $u:M\to N$ is called {\it $F$-harmonic} if it is a critical point of the $F$-energy functional: 
\begin{align}E_F (u) =\int _M F\! \(\frac {|du|^2}{2}\) dv_g.\end{align}  In particular,
if $F (t) =\frac 1p (2t)^{\frac p2}$, then the $F$-energy $E_F(u)$ becomes  $p$-energy, and its critical point $u$ is called {\it $p$-harmonic  map}.
A map $u: (M^{m},g_{M}) \to (N^{n},g_{N})$ is  a \emph{$p$-harmonic morphism} if for any $p$-harmonic function $f$ defined on an open set $V$ of $N^{n}$, the composition $f \circ u$ is $p$-harmonic on $u^{-1}(V)$. 

\subsection{Stress-Energy tensor}\label{S2.3}

Let $(M^{m},g)$ be a smooth Riemannian $m$-manifold. Let $\xi :E\rightarrow M^{m}$ be a smooth Riemannian vector bundle over $(M^{m},g)\, ,$ i.e. a vector bundle such that at each fiber is equipped with a positive inner product $\langle \quad , \quad \rangle_E\, .$
 Set $A^p(\xi )=\Gamma (\Lambda
^pT^{*}M\otimes E)$ the space of smooth $p$-forms on $M^{m}$ with
values in the vector bundle $\xi :E\rightarrow M^{m}$.

For $\omega \in A^p(\xi )$, set $|\omega|^2 = \langle \omega, \omega \rangle$ defined as in (\cite[(2.3)] {DW},). The  authors of \cite{LSC} defined the
following $\mathcal{E}_{F,g}$-energy functional given by
\begin{equation}\notag
\mathcal{E}_{F,g}(\omega )=\int_{M^m}F\(\frac{|\omega |^2}2\)dv_g 
\end{equation}
where $F:[0,+\infty )\rightarrow
[0,+\infty )$ is as before.

The {\it stress-energy associated with the $\mathcal{E}_{F,g}$-energy functional} is
defined as follows:
\begin{equation}\label{2.4}
S_{F,\omega }(X,Y)=F\(\frac{|\omega |^2}2\)g_M(X,Y)-F^{\prime }\(\frac{|\omega |^2%
}2\)u^{\ast}g_N(i_X\omega,i_Y\omega)  
\end{equation}
where $i_X\omega$ is the interior multiplication by the vector field $X$ given by
$$
(i_X\omega)(Y_1,\ldots,Y_{p-1})=\omega(X,Y_1,\ldots,Y_{p-1})
$$
for $\omega \in A^{p}(\xi)$ and any vector fields $Y_{l}$ on $M^{m}$, $1\leq l\leq p-1$.

When $F(t)=t$ and $\omega =du$ for a map $u:M^{m }\rightarrow N^{n}$, $S_{F,\omega }$ is
just the {\it stress-energy tensor} introduced in \cite {BE}.
And when $F(t)=\frac 1n(2t)^{\frac n2}$ and $\omega =du$ for a map
$u:M^{m}\rightarrow N^{n}$, $S_{F,\omega }$ is
 the {\it $n$-stress energy tensor} $S_n$ given by
\begin{equation}\label{2.5}
S_n = \frac {1}{n} |du|^n g_M - |du|^{n-2}u^{\ast} g_N\, .
\end{equation}

\begin{definition} $\omega \in A^p(\xi )\, (p\geq 1)$ is said to satisfy an \emph {$F$-conservation law} if $S_{F,\omega }$ is divergence free, i.e., the $(0,1)$-type tensor field $\operatorname{div} S_{F,\omega }$ vanishes identically {\rm(}i.e.,  $\operatorname{div} S_{F,\omega }\equiv 0${\rm)}.
\end{definition}


The $n$-conservation law is given by \begin{equation}\label{2.6}
\operatorname{div}(S_n) = 0\, \end{equation} (cf. \cite {DW, OC} for details), in which coarea formula was {\it first} employed by Y.\,X. Dong and S.\,W. Wei to derive monotonicity formulas, vanishing theorems, and Liouville theorems on complete noncompact manifolds from conservation laws. 

\section{Main Theorems and Their Proofs}\label{S3}

Assume that $\dim M^{m}=m$ and $\dim N^{n}=n$.

\begin{theorem}\label{T:3.1}
Let $u: (M^{m}, g_M) \to (N^{n}, g_N)$ be a non-constant map. Then the $n$-stress tensor $S_n = 0$ if and only if $m=n$ and $u$ is conformal. 
\end{theorem}
\begin{proof}
If $S_n = 0$, then 
\begin{equation}\label{3.1} u^{\ast}g_{N} = \frac {1}{n} |du|^2 g_{M} = \lambda ^2 g_{M}\end{equation} in the region $du \ne 0$, where $\lambda$ is the dilation and thus
\begin{equation} \begin{aligned} \label{3.2} 0 &= \operatorname{trace}\,  S_n = \frac 1n |du|^n \operatorname{trace}\, g_M - |du|^{n-2} \operatorname{trace}\, u^{*} g_N\\
& = e_n\, m - n e_n\\
& = (m-n)\, e_n,\end{aligned}\end{equation}
where $e_n$ is the $n$-energy density of $u$ given by $e_n = \frac 1n |du|^n$. Hence, we get $m=n$.

Conversely, if $u^{\ast} g_N = \lambda ^2 g_M$ and $m=n$, then we find $$|du|^2 = m \lambda ^2, \, \frac 1n |du|^n = \frac 1n (m \lambda ^2 )^{\frac n2}.$$ Therefore, we obtain
\begin{equation}\label{3.3} S_n =  m^{\frac {n-2}{2}} \frac {(m-n)}{n} {\lambda ^n} g_M=0,
\end{equation} which shows that the $n$-stress tensor $S_n$ vanishes identically.
\end{proof}

\begin{theorem}\label{T:3.2}
If $m > n$ and $u: (M^{m}, g_M) \to (N^{n}, g_N)$ is an $n$-harmonic and conformal map, then  $u$ is homothetic.
\end{theorem}

\begin{proof}
If $u$ is $n$-harmonic, then it follows from \cite [Corollary 2.2]{DW} that $u$ satisfies $n$-conservation law, 
i.e., div$(S_n) = 0$.

In virtue of Theorem \ref{T:3.1} and \eqref{3.3}, with these hypotheses, we find
\begin{equation}\begin{aligned}\label{3.4}
0 = \text{div} (S_n) =  \(m^{\frac {n-2}{2}} \frac{m-n}{n}\)  \text{div} (\lambda^n g_M)
 =  \( m^{\frac{n-2}{2}} \frac{m-n}{n} \) \langle d(\lambda^n), g_M \rangle.
\end{aligned}\end{equation}
Thus, it follows from the assumption $m>n$ that $\lambda$ is a constant. Therefore, $u$  is homothetic.
\end{proof}

Theorems \ref{T:3.2} is sharp in dimensions $m > n$. That is, if $m=n$, then the results no longer hold. Counterexamples can be provided and based on the fact that a conformal map between equal dimensional $n$-manifolds, such as stereographic projections $u: \mathbb E^n \to S^n$ is $n$-harmonic, but $u$ is not homothetic (cf. \cite {WLW, OW}). In fact, Y.\,L. Ou and S.\,W. Wei proved the following:
\vskip.1in

\noindent{\bf Theorem G.} (\cite {OW})
{\it Let $u :  (M^{m},g_{M}) \to (N^{n},g_{N})$ be a non-constant map between Riemannian manifolds with $\dim M = \dim N = n \ge 2$. Then
$u$ is an $n$-harmonic morphism if and only if $u$ is weakly conformal.}
\vskip.1in

While Theorems \ref{T:3.2} on the one hand, augments Theorem G, on the other hand, Theorems \ref{T:3.1} and \ref{T:3.2} generalize
the work of J. Eells and L. Lemaire (\cite {EL2}) in which $n=2$. Furthermore, Theorems \ref{T:3.1} and \ref{T:3.2} augment
a theorem of M. Ara in \cite {A} for the case the zeros of $(n-2) F^{\prime}(t)-2tF^{\prime\prime}(t)$ are being isolated for horizontally conformal $F$-harmonic maps.
Hence we obtain:

\begin{theorem}\label{T:3.3}
 Let $u: (M^{m},g_{M}) \to (N^{n},g_{N})$, $m > n$, be an $F$-harmonic map, which is horizontally conformal  with dilation $\lambda$. 

\noindent
Case 1.  Assume that the zeros of $(n-2) F^{\prime}(t)-2tF^{\prime\prime}(t)$ are isolated. Then the following three properties are equivalent:
\begin{itemize} 
\item[{\rm (1)}] The fibers of $u$ are minimal submanifolds.
 \item[{\rm (2)}]  $\operatorname{grad}(\lambda^2)$ is vertical.
\item[{\rm (3)}] The horizontally distribution of $u$ has mean curvature vector $\frac{\operatorname{grad}(\lambda^2)}{2\lambda^2}$.
\end{itemize}

\noindent 
Case 2. Assume that the zeros of $(n-2) F^{\prime}(t)-2tF^{\prime\prime}(t)$ are not isolated. Then  
\begin{itemize} 
\item[{\rm (1)}] The fibers of $u$ are minimal submanifolds.
 \item[{\rm (2)}] $u$ is homothetic, i.e. $\lambda = C\, ,$ a positive constant.
  \item[{\rm (3)}]  $\operatorname{grad}(\lambda^2) = 0$, hence it is vertical.
 \end{itemize}
\end{theorem}
\begin{proof} Case 1 is exactly   \cite[Theorem 5.1]{A} proved by M. Ara.

For Case 2, statement (1) follows from that fact that general solutions of $$(n-2) F^{\prime}(t)-2tF^{\prime\prime}(t)=0$$ are given by $F(t) = a t^{\frac {n}{2}} +b$ with constants $a,b$. Hence, $u$  is an $n$-harmonic map, and so we may apply Theorem \ref{T:2.1}  to conclude that  fibers of $u$ are minimal in $M^{m}$. Statements
(2) and (3) of Case 2 follow from Theorem \ref{T:3.2} and the fact that $u$ is $n$-harmonic.   
\end{proof}

In examining the converse of Theorem \ref{T:3.3}, Case 2, (1), we characterize the minimal fibers of a horizontally conformal map from the previously untreated case in $F$-harmonic maps:

\begin{theorem}
 Let $u: (M^{m},g_{M}) \to (N^{n},g_{N})$, $m > n$, be a horizontally conformal map. Assume that the zeros of $(n-2) F^{\prime}(t)-2tF^{\prime\prime}(t)$ are not isolated. Then
 the fibers of $u$ are minimal submanifolds if and only if $u$ is an $F$-harmonic map; if and only if $u$ is an $n$-harmonic map. \label{T:3.4}
\end{theorem} 

\begin{proof}
This follows from the fact that when the zeros of $(n-2) F^{\prime}(t)-2tF^{\prime\prime}(t)$ are not isolated, $F$-harmonic map is an $n$-harmonic map, and Theorem \ref{T:2.1}.
\end{proof}

When the target manifold of $u$ is a Riemann surface, i.e. $n=2$, then we associate $u$ with a harmonic map in the following way:

\begin{theorem}
 Let $(N^{2},g_{N})$ be a Riemann surface,  and $u: (M^{m},g_{M}) \to (N^{2},g_{N})$, $m > 2$, be a horizontally conformal map. Assume that the zeros of $-2tF^{\prime\prime}(t)$ are not isolated. Then
 the fibers of $u$ are minimal submanifolds if and only if $u$ is an $F$-harmonic map; if and only if $u$ is a harmonic map. \label{T:3.5}
\end{theorem}

\begin{proof}
This follows from the fact that when the zeros of $-2tF^{\prime\prime}(t)$ are not isolated, $F$-harmonic map is a harmonic map, and Theorem \ref{T:3.4}.
\end{proof}

\section{Applications}\label{S4}

As an application of Theorems \ref{T:3.1} and \ref{T:3.2}, we revisit 

\begin{theorem}[{\bf Theorem C.} (\cite {CW08,CW09}]  \label{T:3.5} Let $u: (M^{m},g_{M}) \to (N^{n},g_{N})$ be an $n$-harmonic morphism which is a submersion. If $N^{n}$ is an orientable manifold and $M^{m}$ is a closed manifold with the $n$-th cohomology class $H^n(M,\mathbf R) = 0$, then the horizontal distribution $\mathcal H$ of $u$ is never integrable. 
\label{T:4.1}\end{theorem}
\begin{proof} Under the hypothesis, in view of Theorem \ref{T:2.1}, $u$ has minimal fibers and, according to Theorem \ref{T:3.2}, 
$\lambda$ is constant.
Let $\{{\bar e}_{1},\ldots,{\bar e}_{n}\}$ be an oriented local orthonormal frame of the base manifold $(N^{n},g_{N})$ and let ${\bar \omega}^{1},\ldots,{\bar \omega}^{n}$ denote the dual 1-forms of $\{{\bar e}_{1},\ldots,{\bar e}_{n}\}$ on $N^{n}$. Then ${\bar \omega}={\bar \omega}^{1}\wedge \cdots  {\bar \omega}^{n}$ is the volume form of $(N^{n},g_{N})$, which is a closed $n$-form on $N^{n}$. 

Consider the pull back of the volume form ${\bar \omega}$ of $N^{n}$ via $u$, which is denoted by  $u^{*}({\bar \omega})$.  Then $u^{*}({\bar \omega})$ is a simple $n$-form on $M^{m}$  satisfying
  \begin{equation}\label{4.1}d \big (u^{*}({\bar \omega})\big )=u^{*}(d{\bar \omega}) =0,\end{equation} due to the fact that the exterior differentiation $d$ and the pullback $u^{*}$ commute. 
  
Assume that $m=\dim M^{m} = n+k$ and let $e_1, \dots, e_{n+k}$ be a local orthonormal frame field with $\omega^1, \dots, \omega^{n+k}$ being its dual coframe fields  on $M^{m}$ such that 
  
\vskip.06in
\noindent   (i) $e_{1}, \dots, e_{n}$ are basic horizontal vector fields satisfying $du(e_{i})=\lambda \bar e_{i}$, $ i=1,\ldots,n$, and  $du(e_{1}), \dots, du(e_{n})$ give a positive orientation of $N^{n}$; and
  
 \noindent  (ii) $e_{n+1}, \dots, e_{n+k}$ are vertical vector fields. 

\vskip.06in
\noindent Then  we have 
\begin{equation}\label{4.2}  \omega^j(e_s) = 0,\;\;\;  \omega^i(e_j)=\delta_{ij},\;\; 1 \le i, j \le n; \; n+1 \le s \le n+k\, . \end{equation} 
Also, it follows from (i) that
\begin{equation}\label{4.3}  u^{\ast} \bar{\omega}^{i}=\frac{1}{\lambda}\omega^{i},\;\; \; i=1,\ldots,n.
 \end{equation} 
  
 If we put 
\begin{equation}\label{4.4} \omega = \omega^1 \wedge \cdots \wedge \omega^{n}\;\;\; {\rm and}\;\;\; \omega^{\bot} = \omega^{n+1}\wedge \cdots \wedge \omega^{n+k}\, , \end{equation}
then 
\begin{equation}\label{4.5}d\omega^{\bot} = \sum _{i=1}^{k} (-1)^{i}\omega^{n+1}\wedge \dots \wedge d\omega^{n+i}\wedge  \cdots \wedge \omega^{n+k}.
\end{equation}
It follows from \e{4.2} and \e{4.5}  that $d\omega^{\bot}=0$ holds identically if and only if the following two conditions are satisfied:
\begin{equation} \label{4.6}d\omega^{\bot}(e_i,e_{n+1}, \dots, e_{n+k}) = 0,\;\;\;  i=1,\ldots, n,\end{equation}
and
\begin{equation} \label{4.7}d\omega^{\bot}(X,Y,V_1,\dots,V_{k-1}) = 0. \end{equation}
for any horizontal vector  fields $X,Y$ and for vertical vector fields $V_{1}, \ldots, V_{k-1}$.

Since the fibers of $u$ are minimal submanifolds of $M^{m}$,  we find  for each $1 \le i \le n$ that
\begin{equation}\begin{aligned}\label{4.8}
d\omega^{\bot}(e_i,\,&e_{n+1}, \dots, e_{n+k})\\
& = \sum _{j=1}^k (-1)^{j+1} \omega^{\bot}([e_i,e_{n+j}], e_{n+1},\dots,\hat{e} _{n+j}, \dots, e_{n+k}) \\
&=\sum _{j=1}^k (-1)^{j+1} \big(\omega^{n+j} (\nabla _{e_i} e_{n+j}) - \omega^{n+j} (\nabla _{e_n+j} e_{i}) \big)\\
&= \sum _{j=1}^k - \langle \nabla _{e_{n+j}} e_{i}, e_{n+j}\rangle\\
&=  \sum _{j=1}^k \langle h(e_{n+j}, e_{n+j}), e_i\rangle \\
&=  0, \end{aligned}
\end{equation}
where ``$\,\hat{\cdot}\,$'' denotes the missing term and $h$ denotes the second fundamental form of fibers in $M$, which prove that condition \e{4.2} holds.

Now, suppose that the horizontal distribution $\mathcal{H}$ is integrable. If $X, Y$ are  horizontal vector fields,   then $[X,Y]$ is also horizontal by Frobenius theorem. So, for vertical vector fields $V_1, \dots, V_k\, ,$ we find (cf. \cite[formula (6.7)]{c2005} or \cite[formula (3.5)] {W9})
\begin{equation} \label{4.9}d\omega^{\bot}(X,Y,V_1,\dots,V_{k-1}) = \omega^{\bot}([X,Y],V_1,\dots,V_{k-1}) = 0. \end{equation}
Consequently, from \e{4.8} and \e{4.9} we get
\begin{equation}\label{4.10}d\omega^{\bot}=0.\end{equation}  

Next, we  show that if  $\mathcal{H}$ is integrable, then 
we have $d\big ((u^{\ast} \bar{\omega}\,)^{\bot}\big )=0\, .$ 
Since $u$ is a horizontally conformal submersion with constant dilation $\lambda$, it preserves orthogonality, which is crucial to horizontal and vertical distributions, and the pullback $u^{\ast}$ expands the length of $1$-form constantly by $\frac 1\lambda$ in every direction. This, via \eqref{4.3} and \eqref{4.10}  leads to
\begin{equation}\begin{aligned}\label{4.11}
d\big ((u^{\ast} \bar{\omega}\,)^{\bot}\big )& 
= d\big ((u^{\ast} \bar{\omega}^{1} \wedge \dots \wedge u^{\ast} \bar{\omega}^{n})^{\bot}\big )  \\
& = d\left(\frac {1}{\lambda} {\omega}^{1} \wedge \cdots \wedge \frac {1}{\lambda} {\omega}^{n}\right)^{\bot}  \\
& = \frac{1}{\lambda^{n}}d\omega^{\bot}
\\& =0.\end{aligned}\end{equation}
Since $d\big ((u^{\ast}\omega)^{\bot}\big )=0$ is equivalent to $u^{\ast}\omega$ being co-closed, it follows  that, under the condition that $\mathcal H$ is integrable, the pullback of the volume form, $u^{\ast}\omega$ is a harmonic $n$-form on $M$. Thus, $u^{\ast}\omega$ gives rise to a non-trivial cohomology class in $H^n(M,\mathbf R)$ by Hodge Theorem (\cite {H}).
Therefore, if $H^n(M,\mathbf R)=0$, then the horizontal distribution $\mathcal{H}$ of $u$ is never integrable.
\end{proof} 

From the proof of Theorem \ref{T:4.1}, we have the following.

\begin{theorem}  Let $u: (M,g_{M}) \to (N,g_{N})$ be an $n$-harmonic morphism  with $ n=\dim  N$ which is a submersion. Then the pull back of the volume element of $N$ is a harmonic $n$-form if and only if the horizontal distribution $\mathcal H$ of $u$ is completely integrable.
\label{T:4.2}\end{theorem}

\end{document}